\documentclass[12pt,letterpaper,twoside]{article}

\usepackage{amsmath,amssymb}
\usepackage{amsthm}
\usepackage{stmaryrd}
\usepackage{fancyhdr}
\usepackage{times}
\usepackage{mathrsfs} % Required for C and H
\usepackage{sectsty}
\usepackage[pdftex,dvips]{geometry}

\sectionfont{\normalsize\bfseries\scshape\centering\MakeUppercase }

\subsectionfont{\normalsize\bfseries\scshape }

\makeatletter
  \def\@seccntformat#1{\csname the#1\endcsname.\quad}
\makeatother

\newtheoremstyle{fact}% name
     {\topsep}%      Space above
     {\topsep}%      Space below
     {\slshape}%         Body font
     {}%         Indent amount (empty = no indent, \parindent = para indent)
     {\bfseries}% Thm head font 
     {}%        Punctuation after thm head
     { }%     Space after thm head: " " = normal interword space;
     {\thmname{#1}\thmnumber{ #2.}\thmnote{ \rm (#3)}}
\newtheorem{theorem}{Theorem}
\newtheorem*{theorem*}{Theorem}

\newtheorem{corollary}[theorem]{Corollary}
\newtheorem{problem}{Problem}
\newtheorem*{problem*}{Problem}

\theoremstyle{definition}

\newtheorem*{remark*}{Remark}

\newtheorem*{question*}{Question}
\newtheorem*{examples*}{Examples}  

\newtheorem*{example*}{Example}

\newtheorem*{convention*}{Convention}

\theoremstyle{fact}

\newtheorem{ftheorem}[theorem]{Theorem}

\def\proofont{\fontseries{bx}\fontshape{sc}\selectfont}
\def\proofname{Proof. }

\newcommand{\Note}[1]{}

\makeatletter
\renewenvironment{proof}[1][\proofname]{\par
  \normalfont
  \topsep6\p@\@plus6\p@ \trivlist
  \item[\hskip\labelsep\noindent\proofont #1]\ignorespaces
}{%
  \qed\endtrivlist
}
\makeatother

\author{G\'abor Luk\'acs
\thanks{I gratefully acknowledge the generous financial support received
from the Killam Trusts and Dalhousie University that enabled me to do 
this research.}}
   
\title{Hereditarily non-topologizable groups
}

\begin{document}

\makeatletter
\let\mytitle\@title
\chead{\small\itshape G. Luk\'acs / \mytitle }
\fancyhead[RO,LE]{\small \thepage}
\makeatother

\maketitle

\def\thanks#1{} 

\thispagestyle{empty}

% References that we may or may not mention explicitly:

\begin{abstract}
A group $G$ is {\em non-topologizable} if the only Hausdorff group 
topology that $G$ admits is the discrete one. Is there an infinite  group 
$G$ such that $H/N$ is non-topologizable for every subgroup $H \leq G$ and  
every normal subgroup $N \vartriangleleft H$? We show that a solution of 
this essentially group theoretic question provides a solution to the 
problem of $c$-compactness.
\end{abstract}

%\paragraph{Introduction.}
Following Ol'shanski\v\i, we say that a group $G$ is {\em 
non-topologizable} if the only Hausdorff group topology that $G$ admits is 
the discrete one. In 1944, Markov asked whether infinite non-topologizable 
groups exist (cf. \cite{Markov}). Markov's problem was solved in 1980, 
independently, by Ol'shanski\v{\i} and Shelah, who constructed infinite 
non-topologizable groups (\cite{Olshanskii} and \cite{Shelah}). Although
Ol'shanski\v\i's example is countable and periodic, Klyachko and Trofimov 
showed that a non-topologizable group need not satisfy either of these 
properties.

\begin{ftheorem}[\cite{KlaTro}] \label{thm:KT}
There exists a torsion-free finitely generated non-topologizable group. 
Thus, there exists a torsion-free non-topologizable group of any 
cardinality.
\end{ftheorem}

(The second statement is obtained from the first one using   
L\"owenheim-Skolem theorem.) In a subsequent paper, Trofimov proved that 
every group embeds into a non-topologizable group of the 
same cardinality (cf.~\cite[Thm.~3]{Trof2}). 

Markov himself obtained a criterion of non-topologizability for countable 
groups with a strong algebraic geometric flavour (Theorem~\ref{thm:cnt:crt} 
below), whose most elegant proof was given by Zelenyuk and Protasov, more 
than half a century later (cf.~\cite{Markov} and~\cite[3.2.4]{PZMono}). 
Given a  monomial 
$f(x)=\nolinebreak g_0x^{k_1}g_2 x^{k_2}g_3\ldots g_{n-1} x^{k_n} g_n$ in 
a single 
variable $x$, with $g_i \in G$ and $k_i \in \mathbb{Z}$, the set $V(f) = 
\{g \in G \mid f(g)=e\}$ is closed in any Hausdorff group 
topology on $G$, because multiplication must be continuous. Thus, if 
$G\backslash \{e\}$ can be represented as $V(f_1) \cup \ldots \cup 
V(f_n)$, where each $f_i$ is a monomial, then $\{e\}$ is open in any 
Hausdorff group topology on $G$, and therefore $G$ is non-topologizable.
In this case, one says that $e$ is {\em algebraically isolated} in $G$. 
The reverse implication holds only if $G$ is countable.

\begin{ftheorem}[\cite{Markov}, {\cite[3.2.4]{PZMono}}] 
\label{thm:cnt:crt}
A countable group $G$ is non-topologizable if and only if $e$ is 
algebraically isolated in $G$.
\end{ftheorem}

Shelah's solution, on the other hand, is uncountable and simple. 
Thus, his result can be rephrased as follows:

\begin{ftheorem}[{\cite{Shelah}}]
Under the Continuum Hypothesis, there is a group $G$ such that 
$G/N$ is non-topologizable for every $N \vartriangleleft G$.
\end{ftheorem}

We say that $G$ is {\em hereditarily non-topologizable} if $H/N$ is 
non-topologizable for every subgroup $H \leq G$ and  every
normal subgroup $N \vartriangleleft H$.
Motivated by Shelah's result, we pose the following problem, and show that 
it is intimately related to the decade-old problem of $c$-compactness of 
topological groups, outlined below.

\begin{problem} \label{prob:GHN}
Is there an infinite hereditarily non-topologizable group?
%{\rm (}resp., countable{\rm )} 
%group $G$ such that 
%$H/N$ is non-topologizable for every subgroup $H \leq G$ and  every 
%normal subgroup $N \vartriangleleft H$?
\end{problem}

%\paragraph{The problem of $\boldsymbol c$-compactness.}
By the well known Kuratowski-Mr\'{o}wka Theorem, a (Hausdorff) topological 
space $X$ is compact if and only if for any (Hausdorff) topological space 
$Y$ the projection $p_Y: X \times Y \rightarrow Y$ is closed. Inspired by 
this theorem, one says that a Hausdorff topological group $G$ is  {\em 
$c$-compact} if for any Hausdorff  group $H$, the image of every  
closed subgroup of $G \times H$ under the projection $\pi_H: G \times H 
\rightarrow H$ is  closed in $H$. The problem of whether 
every $c$-compact topological group is compact has been an open question 
for more than ten years. The most extensive study of $c$-compact 
topological groups was done by Dikranjan and Uspenskij in \cite{DikUsp}, 
which was a source of  inspiration for part of the author's PhD 
dissertation, and his subsequent work (cf.  \cite{GLPHD}, \cite{GL6}, 
\cite{GL4}).

A Hausdorff topological group $G$ is {\em minimal} if there is no coarser 
Hausdorff group topology on $G$ (cf.~\cite{Steph} and~\cite{Dio2}). So, a 
discrete group $G$ is non-topologizable if and only if it is minimal. One 
says that a Hausdorff topological group $G$ is {\em totally minimal} if 
every quotient of $G$ by a closed normal subgroup is minimal 
(cf.~\cite{DikPro}), or equivalently, if every continuous surjective 
homomorphism $f: G \rightarrow H$ is open.  The following two 
results of Dikranjan and Uspenskij link between $c$-compactness and 
total minimality.

\begin{ftheorem}[{\cite[3.6]{DikUsp}}]
Every closed separable subgroup of a $c$-compact group is totally minimal.
\end{ftheorem}

\begin{ftheorem}[{\cite[5.5]{DikUsp}}] \label{thm:DikUsp:discc}
A countable discrete group $G$ is $c$-compact if and only if every 
subgroup of $G$ is totally minimal.
\end{ftheorem}

A discrete group $G$ is hereditarily non-topologizable if and only if the 
discrete topology is totally minimal on every subgroup of $G$. Thus, 
Theorem~\ref{thm:DikUsp:discc} yields: 

\begin{corollary} \label{cor:hnt-cc}
A countable discrete group is $c$-compact if and only if it is 
hereditarily non-topo\-logizable. \qed
\end{corollary}

Recall that a topological group $G$ has {\em small invariant 
neighborhoods} (or briefly, {\em $G$ is SIN}), if any neighborhood $U$ of 
$e \in G$ contains an invariant neighborhood $V$ of $e$, that is, a 
neighborhood $V$ such that $g^{-1} V g = V$ for all $g \in G$. 
Equivalently, $G$ is SIN if its left and right uniformities coincide. 
In a former paper, the author showed that the problem of $c$-compactness 
for locally compact SIN groups can be reduced to the countable discrete 
case (cf.~\cite[4.5]{GL4}). Therefore, Problem~\ref{prob:GHN} is 
equivalent to a special case of the problem of $c$-compactness.

\begin{samepage}
\begin{theorem} \label{thm:red-SIN:cc}
The following statements are equivalent:

\begin{list}{{\rm (\roman{enumi})}}
{\usecounter{enumi}\setlength{\labelwidth}{30pt}\setlength{\topsep}{-8pt}
\setlength{\itemsep}{-4pt}}

\item
every locally compact $c$-compact group admitting small invariant 
neighborhoods is compact;

\item
every countable hereditarily non-topologizable group is finite. \qed

\end{list}
\end{theorem}
\end{samepage}

We conclude with an algebraic consequence of hereditary non-topologizability. 
We denote by  $H^{(k)}$ the $k$-th derived group of a group $H$, that is,
$H^{(1)}=[H,H]$, and $H^{(k)}=[H^{(k-1)},H^{(k-1)}]$.

\begin{theorem} \label{thm:alg}
Let $G$ be a hereditarily non-topologizable group, and let $H \leq G$ be a 
subgroup. Then:

\begin{list}{{\rm (\alph{enumi})}}   
{\usecounter{enumi}\setlength{\labelwidth}{25pt}\setlength{\topsep}{-8pt}
\setlength{\itemsep}{-5pt} \setlength{\leftmargin}{25pt}}

\item
$H$ has a smallest subgroup $N$ of finite index, and $N=[N,N]$;

\item
$H^{(k)}$ has finite index in $H$ for every $k\in \mathbb{N}$;

\item
there is $n \in \mathbb{N}$ such that $H^{(n)}=H^{(n+1)}$;

\item
if $H$ is soluble, then $H$ is finite.

\end{list}

\end{theorem}

In light of Corollary~\ref{cor:hnt-cc}, Dikranjan and Uspenskij's 
results imply Theorem~\ref{thm:alg} (cf.~\cite[3.7-3.12]{DikUsp}). 
Nevertheless, for the sake of completeness, we provide here a direct proof 
that does not rely on the Prodanov-Stoyanov theorem (cf.~\cite{ProdStoj}).

\begin{proof}
(a) Let $\{N_\alpha\}$ be the collection of normal subgroups of finite 
index in $H$, and set\linebreak[4] $N=\bigcap N_\alpha$.  Since $G$ is 
hereditarily 
non-topologizable, the discrete topology is the only Hausdorff group 
topology on $H/N$. On the other hand, $H/N$ embeds into the product 
$P=\prod H/N_\alpha$, and $P$ is compact Hausdorff, because 
each $H/N_\alpha$ is finite. Thus, the image of $H/N$ in $P$ can be 
discrete only if it is finite. Therefore, $N$ has finite index in $H$.
(The closure of the image of $H/N$ in $P$ is called the {\em pro-finite 
completion} of $H/N$.)

Let $H_1 \leq H$ be such that $|H:H_1|=l$. Then $H_1$ contains a normal 
subgroup $N_1$ of $H$ such that $|H:N_1|\leq l!$, and hence
$N \leq N_1 \leq H_1$, as desired. By (a), $[N,N]$ has finite index 
in $N$, and so the last statement follows.

(b) Let $A=H/[H,H]$ be the maximal abelian quotient of $H$. We denote by 
$\hat A$ the 
group $\hom_\mathbb{Z}(A,\mathbb{T})$, where 
$\mathbb{T}=\mathbb{R}/\mathbb{Z}$.
Consider the homomorphism 
\begin{align}
\rho_A\colon A & \longrightarrow \mathbb{T}^{\hat A} \\
a & \longmapsto (\chi(a))_{\chi \in \hat A}.
\end{align}
For every non-zero $a \in A$,  there is a homomorphism 
$\chi_a\colon \langle a \rangle \rightarrow \mathbb{T}$ such that 
$\chi_a(a)\neq 0$. Since $\mathbb{T}$ is injective, 
$\chi_a$ can be extended to $\bar\chi_a \colon A \rightarrow \mathbb{T}$.
Thus, $\rho_A(a)\neq 0$, and so $\rho_A$ is injective. The {\em Bohr 
topology} on $A$ is the initial topology induced by $\rho_A$, that is, the 
subgroup topology when $A$ is viewed as a subgroup of 
$\mathbb{T}^{\hat A}$. Since $G$ is hereditarily non-topologizable, the 
only Hausdorff group topology on $A$ is the discrete one. On the other hand, 
the subgroup $\rho_A(A)$ of the compact Hausdorff group $\mathbb{T}^{\hat 
A}$ is discrete if and only if $A$ is finite. Therefore, $A$ is finite. 
Hence, the statement follows by an inductive reiteration of this argument.

(c) and (d) follow from (a) and (b).
\end{proof}

\section*{Acknowledgments}

My deepest thanks go to Walter Tholen, for introducing me to the problem 
of $c$-compactness, and for his attention to my work.

\medskip\noindent
I wish to thank Robert Par\'e for the valuable discussions that were of 
great assistance in writing this paper.

{\footnotesize

\bibliography{notes,notes2,notes3}

}

\begin{samepage}

{\bigskip\bigskip\noindent 
Department of Mathematics and Statistics\\
Dalhousie University\\
Halifax, B3H 3J5, Nova Scotia\\
Canada

\nopagebreak
\bigskip\noindent{\em e-mail: lukacs@mathstat.dal.ca} }
\end{samepage}

\end{document}